\documentclass[11pt]{article}
\usepackage[a4paper,margin=1.4in]{geometry}

\usepackage[colorlinks,citecolor=magenta,linkcolor=black]{hyperref}
\pdfpagewidth=\paperwidth \pdfpageheight=\paperheight
\usepackage{amsfonts,amssymb,amsthm,amsmath,eucal,tabu,url}
\usepackage{pgf}
 \usepackage{array}
 \usepackage{pstricks}
 \usepackage{pstricks-add}
 \usepackage{pgf,tikz}
 \usetikzlibrary{automata}
 \usetikzlibrary{arrows}
 \usepackage{indentfirst}
\usepackage{tabularx} 
\usepackage{longtable}

\theoremstyle{plain}
\newtheorem{thm}{Theorem}[section]
\newtheorem{theorem}[thm]{Theorem}

\newtheorem{conjecture}[thm]{Conjecture}
\newtheorem*{theoremA}{Main Theorem}

\theoremstyle{definition}
\newtheorem{definition}[thm]{Definition}
\newtheorem{remark}[thm]{Remark}

\newtheorem{thevarthm}[thm]{\varthmname}

\newenvironment{varthm*}[1]{\trivlist\item[]{\bf #1.}\it}{\endtrivlist}

\renewcommand\geq{\geqslant}

\renewcommand\leq{\leqslant}

\newcommand\be{\begin{eqnarray*}}
\newcommand\ee{\end{eqnarray*}}

\newcommand\cala{{\mathcal A}}

\newcommand\newop[2]{\def#1{\mathop{\rm #2}\nolimits}}
\newop\Der{Der}
\newop\pdeg{pdeg}
\newop\Jac{Jac}
\newop\log{log}
\newop\ord{ord}
\newop\Gal{Gal}
\newop\SL{SL}
\newop\Bl{Bl}
\newop\mult{mult}
\newop\mass{mass}
\newop\div{div}
\newop\codim{codim}
\newop\sing{sing}
\newop\vdim{vdim}
\newop\edim{edim}
\newop\Ass{Ass}
\newop\size{size}
\newop\reg{reg}
\newop\satdeg{satdeg}
\newop\supp{supp}
\newop\Neg{Neg}
\newop\Nef{Nef}
\newop\Nefh{Nef_H}
\newop\Eff{Eff}
\newop\Zar{Zar}
\newop\MB{MB}
\newop\MBxC{MB\mathit{(x,C)}}
\newop\NnB{NnB}
\newop\Bigg{Big}
\newop\Effbar{\overline{\Eff}}

\def\keywordname{{\bfseries Keywords}}%
\def\keywords#1{\par\addvspace\medskipamount{\rightskip=0pt plus1cm
\def\and{\ifhmode\unskip\nobreak\fi\ $\cdot$
}\noindent\keywordname\enspace\ignorespaces#1\par}}

\newcommand\K{\mathbb K}

\def\subclassname{{\bfseries Mathematics Subject Classification
(2020)}\enspace}
\def\subclass#1{\par\addvspace\medskipamount{\rightskip=0pt plus1cm
\def\and{\ifhmode\unskip\nobreak\fi\ $\cdot$
}\noindent\subclassname\ignorespaces#1\par}}

\begin{document}
\title{On the nearly free simplicial line arrangements with up to $27$ lines}
\author{Marek Janasz}
\date{\today}
\maketitle
\thispagestyle{empty}
\begin{abstract}
 In the present note we provide a complete classification of nearly free (and not free simultaneously) simplicial arrangements of $d\leq 27$ lines.
\subclass{14N20, 14C20}
\end{abstract}
\section{Introduction}
The theory of line arrangements is a classical subject of studies in many branches of contemporary mathematics. In the recent years, many authors wanted to understand possible linkages between combinatorial and geometric properties of line arrangements. Let us recall that the famous Terao's conjecture predicts that the so-called freeness of a given arrangement of lines $\mathcal{A}$ is determined by the intersection poset of $\mathcal{A}$. It is very difficult to predict whether Terao's conjecture is true, and in order to approach this problem Dimca and Sticlaru in \cite{DimcaSticlaru} defined a new class curves which is called nearly free. This class is designed as a natural generalization of free curves and it is important in the context of a potential counterexample to Terao's conjecture. It seems that the class of nearly free arrangements is more accessible, and it is definitely much wider. In the present note, which can be considered as an appendix to works devoted to simplicial line arrangements in the real projective plane, we want to understand which sporadic examples of simplicial line arrangements in the real projective plane are nearly free and not free. Even if the classification problem of simplicial line arrangements is open in its whole generality, we will use a great result due to M. Cuntz which provides a complete classification of simplicial arrangements up to $27$ lines and, in this way, we provide a complete classification result of nearly free sporadic simplicial arrangements up to $27$ lines. Our main result, surprising to us, can be formulated as follows.

\begin{theoremA}
A sporadic simplicial line arrangement $\mathcal{A} \subset \mathbb{P}^{2}_{\mathbb{R}}$ is nearly free if and only if $\mathcal{A} = \mathcal{A}(17,6)$ according to Cuntz's catalogue.
\end{theoremA}

\begin{remark}
More precisely, $\mathcal{A}(17,6)$ is a sporadic simplicial line arrangement consisting of $17$ lines and it has $16$ double, $15$ triple, $10$ quadruple, and one sixtuple intersection point.
\end{remark}

It means that the class of free sporadic simplicial line arrangements is barely different from the class of nearly free sporadic simplicial line arrangements provided that we restrict our attention to $d\leq 27$ lines.

In order to prove Main Theorem, we will use combinatorial properties of the singular points of sporadic simplicial line arrangements. This allows us to determine all those sporadic arrangements for which the total Milnor number is determined exclusively by a polynomial equation of degree $2$ that depends only on the number of lines and the minimal degree of the syzygies between partial derivatives of the defining polynomial. In the last step, using cohomological methods, we are able to determine those arrangements which are purely nearly free.

The structure of the paper goes as follows. In Section 2, we provide all necessary definitions and tools related to simplicial and nearly free line arrangements. In Section 3, we provide our proof of Main Theorem. All necessary symbolic computations were performed with use of \verb{Singular{ \cite{Singular}.
\section{Preliminaries}
In the section, we recall all necessary notations and definitions. For more information in this area please consult \cite{Dimca,OT92}.

Let $\K$ be any field and consider $S:=\K[x,y,z]$ the graded polynomial ring over $\mathbb{K}$.
\begin{definition}
A finite collection of $d$ lines $\mathcal{L} = \{\ell_{1}, ..., \ell_{d} \}\subset \mathbb{P}^{2}_{\mathbb{K}}$ is called an arrangement of lines in the projective plane over $\mathbb{K}$.
\end{definition}
For an arrangement $\mathcal{L} = \{\ell_{1}, ..., \ell_{d}\}$ we denote by ${\rm Sing}(\mathcal{L})$ the set of all intersection points among the lines, i.e., points in the plane where at least two lines from $\mathcal{L}$ meet, and for such an intersection point $p \in {\rm Sing}(\mathcal{L})$ we denote by ${\rm mult}_{p}$ its multiplicity, i.e., the number of lines passing through the point $p$. Following Hirzebruch's convention, we denote by $t_{r}$ the number of all intersection points of multiplicity $r\geq 2$. 

We define the class of simplicial line arrangements in the real projective plane via Melchior's result \cite{Melchior}.

\begin{definition}
Let $\mathcal{L} = \{\ell_{1}, ..., \ell_{d}\} \subset \mathbb{P}^{2}_{\mathbb{R}}$ of $d\geq 3$ lines such that $t_{d} = 0$. Then $\mathcal{L}$ is a simplicial line arrangement if and only if
$$t_{2} = 3 + \sum_{r\geq 4}(r-3)t_{r}.$$
\end{definition}
Classically, a simplicial line arrangement $\mathcal{L} \subset \mathbb{P}^{2}_{\mathbb{R}}$ is an arrangement for which all connected components of the complement $M(\mathcal{L}) := \mathbb{P}^{2}_{\mathbb{R}} \setminus \mathcal{L}$ are \textbf{open triangles}. It is worth recalling that simplicial line arrangements were studied, for may years, by Gr\"unbaum, and he discovered three infinite families of such arrangements and around $90$ additional examples which are nowadays called \textbf{sporadic}. The collection of the three infinite families and around $90$ sporadic examples is called in the literature as Gr\"unbaum's catalogue. One of the most important conjectures related to simplicial line arrangements is motivated by a strong claim of Gr\"unbaum \cite[p. 4]{Gru}.
\begin{conjecture}
Except only finitely many corrections, Gr\"unbaum's catalogue is complete.
\end{conjecture}
In other words, one expects that there are only three infinite families of simplicial line arrangements. A stronger conjecture, proposed by Cuntz and Geis in \cite[Conjecture 1.6]{CuntzGeis}, predicts even more.

\begin{conjecture}
Let $\mathcal{L}$ be a sporadic simplicial line arrangement in $\mathbb{P}^{2}_{\mathbb{R}}$ of $d$ lines. Then $d\leq 37$.
\end{conjecture}

The main aim of the present note is to understand the homological properties of Jacobian ideals given by simplicial line arrangements.  In order to do so, let recall some crucial definitions. For a reduced curve $C \subset \mathbb{P}^{2}_{\mathbb{C}}$ of degree $d$ given by $f = 0$ we denote by $J_{f} = \langle \partial_{x}\,  f, \, \partial_{y} \, f,\partial_{z} \, f \rangle$ the Jacobian ideal and by $\mathfrak{m} = \langle x,y,z \rangle$ the irrelevant ideal. Consider the graded $S$-module $N(f) = I_{f} / J_{f}$, where $I_{f}$ is the saturation of $J_{f}$ with respect to $\mathfrak{m}$.

\begin{definition}
We say that a reduced plane curve $C$ is \emph{nearly free} if $N(f) \neq 0$ and for every $k$ one has ${\rm dim} \, N(f)_{k} \leq 1$. 
\end{definition}

Recall that  for a curve $C$ given by $f \in S$ we define the Milnor algebra as $M(f) = S / J_{f}$.
The description of $M(f)$ for nearly free curves comes from \cite{DimcaSticlaru} as follows.

\begin{theorem}[Dimca-Sticlaru]
\label{DimSti}
If $C$ is a nearly free curve of degree $d$ given by $f \in S$, then the minimal free resolution of the Milnor algebra $M(f)$ has the following form:
\begin{equation*}
\begin{split}
0 \rightarrow S(-b-2(d-1))\rightarrow S(-d_{1}-(d-1))\oplus S^2(-d_{2}-&(d-1))  \\& \rightarrow S^{3}(-d+1)\rightarrow S
\end{split}
\end{equation*} for some integers $d_{1},d_{2}, b$ such that $d_{1} + d_{2} = d$ and $b=d_{2}-d+2$. In that case, the pair $(d_{1},d_{2})$ is called the set of exponents of $C$.
\end{theorem}


The nearly freeness can be also studied via the following result due to Dimca \cite[Theorem 1.3]{Dimca1}, and this result is a vital technical tool for our proposes.
\begin{theorem}[Dimca]
Let $\mathcal{L}\subset \mathbb{P}^{2}_{\mathbb{C}}$ be an arrangement of $d$ lines and let $f=0$ be its defining equation. Denote by $r$ the minimal degree among all the Jacobian relations, i.e., the minimal degree $r$ for the triple $(a,b,c) \in S_{r}^{3}$ such that $a \cdot \partial_{x}(f) + b \cdot \partial_{y}(f) + c\cdot \partial_{z}(f) = 0$. Assume that $r\leq d/2$, then $\mathcal{L}$ is nearly free if and only if
\begin{equation}
\label{Milnor}
r^2 - r(d-1) + (d-1)^2 = \mu(\mathcal{L})+1,
\end{equation}
where $\mu (\mathcal{C})$ is the total Milnor number of $\mathcal{L}$, i.e.,
$$\mu(\mathcal{L}) = \sum_{p \in {\rm Sing}(\mathcal{L})} ({\rm mult}_{p}-1)^{2}.$$
\end{theorem}
Finally, let us also present a cohomological description of free arrangements, see \cite{DimcaSticlaru} for details.
\begin{theorem}
Let $C \subset \mathbb{P}^{2}_{\mathbb{C}}$ be a reduced curve of degree $d$ and let $f=0$ be its defining equation. Then $C$ is free if and only if then the minimal free resolution of the Milnor algebra $M(f)$ has the following form:
\begin{equation*}
0 \rightarrow S(-d_{1}-(d-1)) \oplus S(-d_{2}-(d-1)) \rightarrow S^{3}(-d+1)\rightarrow S
\end{equation*}
with $d_{1} + d_{2} = d-1$.  The pair $(d_{1},d_{2})$ is called the set of exponents of $C$.
\end{theorem}
\section{Proof of Main Result}
\begin{proof}
Here we want to present the main idea standing behind our proof. First of all, the table below presents all known sporadic simplicial line arrangements in the real projective plane having at most $27$ lines. We have, according to Cuntz's catalogue, around $70$ such arrangements. In the table below we provide additionally the total Milnor number of a given arrangement $\mathcal{A}(x,y)$ (here $x$ denotes the number of lines in the  given arrangement and $y$ its type), the discriminant $\triangle_{r}$ for  (\ref{Milnor}) computed with respect to $r$ as variable, and we provide information about the roots of (\ref{Milnor}) computed with respect to $r$.

Here is the outline of our strategy:
\begin{itemize}
    \item Among all sporadic simplicial line arrangements we detect those for which $\sqrt{\triangle_{r}}$ is an integer.
    \item For those line arrangements with an integral value of $\sqrt{\triangle_{r}}$, we extract all arrangements for which (\ref{Milnor}), computed with respect to $r$, has integral roots.
    \item Finally, after the above two-step process, we compute the minimal free resolutions of Milnor algebras, minimal degrees of the Jacobian relations and, based on that information, we detect those sporadic arrangements which are nearly free and not free.
\end{itemize}

We start with the aforementioned table.
{\footnotesize
\begin{center}
\begin{longtable}{|l|l|l|l|l|}
\caption{The list of sporadic simplicial line arrangements up to $27$ lines.} \label{tab:long} \\

\hline \multicolumn{1}{|c|}{\textbf{$\cala(n,k)$}} & \multicolumn{1}{c|}{$(t_2,t_3,\ldots)$} & \multicolumn{1}{c|}{$\mu(\mathcal{L})$} & \multicolumn{1}{c|}{$\triangle_{r}$} & \multicolumn{1}{c|}{roots}\\ \hline 
\endfirsthead
 
\multicolumn{4}{c}%
{{\bfseries \tablename\ \thetable{} -- continued from the previous page }} \\
\hline \multicolumn{1}{|c|}{\textbf{$\cala(n,k)$}} & \multicolumn{1}{c|}{$(t_2,t_3,\ldots)$} & \multicolumn{1}{c|}{$\mu(\mathcal{L})$} & \multicolumn{1}{c|}{$\triangle_{r}$}& \multicolumn{1}{c|}{roots}\\ \hline 
\endhead

\hline
\endfoot

\hline \hline
\endlastfoot
$\cala(7,1)$ & $(3,6)$ & 27 & $>0$ & $r_1=0,\:r_2=2$ \\
$\cala(9,1)$ & $(6,4,3)$ & 49 & $>0$ &  real\\
$\cala(10,2)$ & $(6,7,3)$ & 61 & $>0$ &  real\\
$\cala(10,3)$ & $(6,7,3)$ & 61 & $>0$ &  real\\
$\cala(11,1)$ & $(7,8,4)$ & 75 & $>0$ & $r_1=4,\:r_2=6$\\
$\cala(12,2)$ & $(8,10,3,1)$ & 91 & $>0$ &  real\\
$\cala(12,3)$ & $(9,7,6)$ & 91 & $>0$ &  real\\
$\cala(13,2)$ & $(12,4,9)$ & 109 & $>0$ &  real\\
$\cala(13,3)$ & $(10,10,3,2)$ & 109 & $>0$ &  real\\
$\cala(13,4)$ & $(6,18,3)$ & 105 & $<0$ &  complex\\
$\cala(14,2)$ & $(11,12,4,2)$ & 127 & $>0$ &  real\\
$\cala(14,3)$ & $(9,16,4,1)$ & 125 & $<0$ &  complex\\
$\cala(14,4)$ & $(10,14,4,0,1)$ & 127 & $>0$ &  real\\
$\cala(15,1)$ & $(15,10,0,6)$ & 151 & $>0$ &  real\\
$\cala(15,2)$ & $(13,12,6,2)$ & 147 & $>0$ & $r_1=6,\:r_2=8$\\
$\cala(15,3)$ & $(12,13,9)$ & 145 & $<0$ &  complex\\
$\cala(15,4)$ & $(12,14,6,0,1)$ & 147 & $>0$ & $r_1=6,\:r_2=8$\\
$\cala(15,5)$ & $(9,22,0,3)$ & 145 & $<0$ &  complex\\
$\cala(16,2)$ & $(14,15,6,1,1)$ & 169 & $>0$ &  real\\
$\cala(16,3)$ & $(15,13,6,3)$ & 169 & $>0$ &  real\\
$\cala(16,4)$ & $(15,15,0,6)$ & 171 & $>0$ &  real\\
$\cala(16,5)$ & $(14,16,3,4)$ & 169 & $>0$ &  real\\
$\cala(16,6)$ & $(15,12,9,0,1)$ & 169 & $>0$ &  real\\
$\cala(16,7)$ & $(12,19,6,0,1)$ & 167 & $<0$ &  complex\\
$\cala(17,2)$ & $(16,16,7,0,2)$ & 193 & $>0$ &  real\\
$\cala(17,3)$ & $(18,12,7,4)$ & 193 & $>0$ &  real\\
$\cala(17,4)$ & $(16,16,7,0,2)$ & 193 & $>0$ &  real\\
$\cala(17,5)$ & $(16,18,1,6)$ & 193 & $>0$ &  real\\
$\cala(17,6)$ & $(16,15,10,0,1)$ & 191 & $0$ & $r_0=8$\\
$\cala(17,7)$ & $(13,22,7,0,1)$ & 189 & $<0$ &  complex\\
$\cala(17,8)$ & $(14,20,7,2)$ & 189 & $<0$ &  complex\\
$\cala(18,2)$ & $(18,18,6,3,1)$ & 217 & $>0$ &  real\\
$\cala(18,3)$ & $(19,16,6,5)$ & 217 & $>0$ &  real\\
$\cala(18,4)$ & $(18,19,3,6)$ & 217 & $>0$ &  real\\
$\cala(18,5)$ & $(18,19,3,6)$ & 217 & $>0$ &  real\\
$\cala(18,6)$ & $(18,16,12,0,1)$ & 215 & $<0$ &  complex\\
$\cala(18,7)$ & $(18,18,6,3,1)$ & 217 & $>0$ &  real\\
$\cala(18,8)$ & $(16,22,6,2,1)$ & 215 & $<0$ &  complex\\
$\cala(19,1)$ & $(21,18,6,0,4)$ & 247 & $>0$ &  real\\
$\cala(19,2)$ & $(21,18,6,6)$ & 243 & $>0$ & $r_1=8,\:r_2=10$\\
$\cala(19,3)$ & $(24,12,6,6,1)$ & 247 & $>0$ &  real\\
$\cala(19,4)$ & $(20,20,6,4,1)$ & 243 & $>0$ & $r_1=8,\:r_2=10$\\
$\cala(19,5)$ & $(20,20,6,4,1)$ & 243 & $>0$ & $r_1=8,\:r_2=10$\\
$\cala(19,6)$ & $(20,20,6,4,1)$ & 243 & $>0$ & $r_1=8,\:r_2=10$\\
$\cala(19,7)$ & $(21,15,15,0,1)$ & 241 & $<0$ &  complex\\
$\cala(20,2)$ & $(25,15,10,6)$ & 271 & $>0$ &  real\\
$\cala(20,3)$ & $(21,24,6,4,0,1)$ & 271 & $>0$ &  real\\
$\cala(20,4)$ & $(23,20,7,5,1)$ & 271 & $>0$ &  real\\
$\cala(20,5)$ & $(20,26,4,4,0,0,1)$ & 273 & $>0$ &  real\\
$\cala(21,2)$ & $(30,10,15,6)$ & 301 & $>0$ &  real\\
$\cala(21,3)$ & $(24,24,9,0,4)$ & 301 & $>0$ &  real\\
$\cala(21,4)$ & $(22,28,6,4,0,0,1)$ & 301 & $>0$ &  real\\
$\cala(21,5)$ & $(26,20,9,4,2)$ & 301 & $>0$ &  real\\
$\cala(21,6)$ & $(25,20,15,2,1)$ & 297 & $<0$ &  complex\\
$\cala(21,7)$ & $(24,22,15,3)$ & 295 & $<0$ &  complex\\
$\cala(22,2)$ & $(24,30,12,3,1)$ & 325 & $<0$ &  complex\\
$\cala(22,3)$ & $(27,28,0,12)$ & 331 & $>0$ &  real\\
$\cala(22,4)$ & $(27,25,9,3,3)$ & 331 & $>0$ &  real\\
$\cala(22,5)$ & $(12,58,0,0,3)$ & 319 & $<0$ &  complex\\
$\cala(23,1)$ & $(27,32,10,4,2)$ & 359 & $<0$ &  complex\\
$\cala(23,2)$ & $(16,56,2,0,1,2)$ & 355 & $<0$ &  complex\\
$\cala(24,2)$ & $(32,32,0,12,0,0,1)$ & 401 & $>0$ &  real\\
$\cala(24,3)$ & $(31,32,9,5,3)$ & 395 & $<0$ &  complex\\
$\cala(24,4)$ & $(20,54,4,0,0,2,1)$ & 393 & $<0$ &  complex\\
$\cala(25,2)$ & $(36,28,15,0,6)$ & 433 & $>0$ &  real\\
$\cala(25,3)$ & $(30,40,15,6)$ & 421 & $<0$ &  complex\\
$\cala(25,4)$ & $(36,30,9,6,4)$ & 433 & $>0$ &  real\\
$\cala(25,5)$ & $(36,32,0,8,4,0,1)$ & 441 & $>0$ &  real\\
$\cala(25,6)$ & $(36,30,9,6,4)$ & 433 & $>0$ &  real\\
$\cala(25,7)$ & $(33,34,12,2,3,0,1)$ & 433 & $>0$ &  real\\
$\cala(25,8)$ & $(24,52,6,0,0,0,3)$ & 433 & $>0$ &  real\\
$\cala(26,2)$ & $(35,40,10,11)$ & 461 & $<0$ &  complex\\
$\cala(26,3)$ & $(37,36,9,6,3,1)$ & 469 & $>0$ &  real\\
$\cala(26,4)$ & $(35,39,10,4,3,0,1)$ & 469 & $>0$ &  real\\
$\cala(27,1)$ & $(40,40,6,14,1)$ & 503 & $<0$ &  complex\\
$\cala(27,2)$ & $(39,40,10,6,2,2)$ & 507 & $>0$ & $r_1=12,\:r_2=14$\\
$\cala(27,3)$ & $(39,40,10,6,2,2)$ & 507 & $>0$ & $r_1=12,\:r_2=14$\\
$\cala(27,4)$ & $(38,42,9,6,3,0,1)$ & 507 & $>0$ & $r_1=12,\:r_2=14$\\
\end{longtable}
\end{center}
}

Based on what we have seen so far, we can check directly that the following arrangements pass the first two steps of our selection, namely:
$$\mathcal{A}(7,1), \mathcal{A}(11,1),\mathcal{A}(15,2), \mathcal{A}(15,4), \mathcal{A}(17,6), \mathcal{A}(19,2), \mathcal{A}(19,4), \mathcal{A}(19.5), \mathcal{A}(19,6),$$ $$ \mathcal{A}(27,2), \mathcal{A}(27,3),\mathcal{A}(27,4).$$
Now, according to Step 3, we present a detailed discussion regarding nearly freeness and freeness of the extracted arrangements.
\begin{itemize}
\item[$\mathcal{A}(7,1): $] The minimal free resolution of the Milnor algebra has the following form
    $$0\rightarrow S^{2}(-9) \rightarrow S^{3}(-6) \rightarrow S ,$$
 which means that $\mathcal{A}(7,1)$ is free.
\item[$\mathcal{A}(11,1): $] The minimal free resolution of the Milnor algebra has the following form
    $$0\rightarrow S^{2}(-15) \rightarrow S^{3}(-10) \rightarrow S ,$$
so $\mathcal{A}(11,1)$ is free.
\item[$\mathcal{A}(15,2): $] The minimal free resolution of the Milnor algebra has the following form
    $$0\rightarrow S^{2}(-21) \rightarrow S^{3}(-14) \rightarrow S ,$$
so $\mathcal{A}(15,2)$ is free.
\item[$\mathcal{A}(15,4): $] The minimal free resolution of the Milnor algebra has the following form
$$0 \rightarrow S^{2}(-21) \rightarrow S^{3}(-14) \rightarrow S.$$
so $\mathcal{A}(15,4)$ is free.
\item[$\mathcal{A}(17,6): $] The minimal free resolution of the Milnor algebra has the following form
     $$0\rightarrow S(-26) \rightarrow S^{2}(-25) \oplus S(-24) \rightarrow S^{3}(-16) \rightarrow S.$$
Since the minimal degree of the Jacobian relations $r$ is equal to $8$ and it satisfies Equation (\ref{Milnor}), \textbf{then $\mathcal{A}(17,6)$ is nearly free}.
\item[$\mathcal{A}(19,2): $] The minimal free resolution of the Milnor algebra has the following form
    $$0\rightarrow S^{2}(-27) \rightarrow S^{3}(-18) \rightarrow S ,$$
so $\mathcal{A}(19,2)$ is free.
\item[$\mathcal{A}(19,4): $] The minimal free resolution of the Milnor algebra has the following form
    $$0\rightarrow S^{2}(-27) \rightarrow S^{3}(-18) \rightarrow S ,$$
so $\mathcal{A}(19,4)$ is free.
\item[$\mathcal{A}(19,5): $] The minimal free resolution of the Milnor algebra has the following form
    $$0\rightarrow S^{2}(-27) \rightarrow S^{3}(-18) \rightarrow S ,$$
so  $\mathcal{A}(19,5)$ is free.
\item[$\mathcal{A}(19,6): $] The minimal free resolution of the Milnor algebra has the following form
    $$0\rightarrow S^{2}(-27) \rightarrow S^{3}(-18) \rightarrow S ,$$
so $\mathcal{A}(19,6)$ is free.
\item[$\mathcal{A}(27,2): $] The minimal free resolution of the Milnor algebra has the following form
    $$0\rightarrow S(-51) \rightarrow S(-49) \oplus S(-41) \oplus S(-39) \rightarrow S^{3}(-26) \rightarrow S ,$$
so according to Theorem \ref{DimSti} arrangement $\mathcal{A}(27,2)$ is not nearly free.
\item[$\mathcal{A}(27,3): $] The minimal free resolution of the Milnor algebra has the following form
    $$0\rightarrow S(-51) \rightarrow S(-49) \oplus S(-41) \oplus S(-39) \rightarrow S^{3}(-26) \rightarrow S ,$$
so according to Theorem \ref{DimSti} arrangement $\mathcal{A}(27,3)$ is not nearly free.
\item[$\mathcal{A}(27,4): $] The minimal free resolution of the Milnor algebra has the following form
    $$0\rightarrow S^{2}(-39) \rightarrow S^{3}(-26) \rightarrow S ,$$
    so $\mathcal{A}(27,4)$ is free.
\end{itemize}
This completes the proof.
\end{proof}
\section*{Acknowledgments}
I would like to thank Piotr Pokora for his guidance during the project and useful suggestions.

\vskip 0.5 cm

\bigskip
Marek Janasz,
Department of Mathematics,
Pedagogical University of Krakow,
ul. Podchorazych 2,
PL-30-084 Krak\'ow, Poland. \\
\nopagebreak
\textit{E-mail address:} \texttt{marek.janasz@up.krakow.pl}

\end{document}